\input amsppt.sty
\loadbold
\TagsOnRight
\hsize 30pc
\vsize 47pc
\magnification 1200

\redefine\g{\frak g}

\redefine\a{\frak a}
\redefine\B{{\Cal B}}
\redefine\D{{\Cal D}}

\redefine\l{\frak l}
\redefine\k{\frak k}
\redefine\h{\frak h}
\redefine\m{\frak m}
\redefine\n{\frak n}
\redefine\p{\frak p}
\redefine\t{\frak t}
\redefine\z{\frak z}
\redefine\C{\Bbb C}
\redefine\R{\Bbb R}

\redefine\Z{Z_\omega}

\redefine\1{{1\over 2}}
\redefine\der{{d\over{dt}}}

\redefine\d{{\text{d}}}

\redefine \bJ{{\Bbb J}}

\define\={\overset\text{def}\to=}

\define\bV{\Bbb V}

\topmatter
\title K\"ahler manifolds with large isometry group\endtitle
\author  Fabio Podest\`a and Andrea Spiro \endauthor 
\address Math. Institute, Via dell'Agnolo, I- Firenze\quad -\quad 
Math. Department, Via delle Brecce Bianche, I- Ancona\endaddress
\email podesta\@ udini.math.unifi.it\quad spiro\@ anvax1.cineca.it
\endemail 
\keywords K\"ahler-Einstein metrics, cohomogeneity one actions\endkeywords
\subjclass  53C55, 53C25, 57S15\endsubjclass
\thanks \endthanks 
 \abstract We investigate compact K\"ahler manifolds, which 
are acted on by a semisimple 
compact Lie group $G$ of 
isometries with one hypersurface orbit. In case of ordinary action and
projectable complex structure, we set up a one to 
one correspondence between such manifolds and
abstract models. The Ricci tensor is then computed and 
we fully characterize the K\"ahler-Einstein
manifolds with an ordinary cohomogeneity one action and projectable 
complex structure. 
\endabstract

\leftheadtext{\smc K\"ahler manifolds with large isometry group}
\rightheadtext{\smc F. Podest\`a - A. Spiro}

\endtopmatter
\document
\subhead 1. Introduction \endsubhead
\bigskip
A Riemannian manifold is called  of cohomogeneity one 
if it is acted on by a closed Lie group $G$ of isometries
with principal orbits of codimension one. This class of manifolds
have at least two good reasons for being considered
particularly appealing: their degree of symmetry is 
so high that classification theorems of purely algebraic 
nature are still possible in several situations (see e.g. [Br],
[AA], [AA1], [HL]); at the same time, they allow to construct 
non homogeneous examples of Riemannian manifolds with 
special geometric properties, like Einstein metrics, exceptional holonomy (see
e.g. [Be]).\par
We are interested in compact non homogeneous
K\"ahler-Einstein manifolds and several of them have been 
already constructed by Koiso and Sakane in the cohomogeneity
one category (see e.g. [Sa], [KS]).  The aims of the
present paper consist in giving an explicit classification
of compact cohomogeneity one K\"ahler manifolds with 
vanishing first Betti number and in using it for obtaining
a complete list of the K\"ahler-Einstein manifolds in that 
family.\par
It is well known (see e.g. [GS]) 
that  the vanishing of $b_1(M)$ implies the existence of a 
$G$-equivariant moment
mapping $\mu \: M \to \g^*$ and this fact has an important
consequence on the algebraic structure of $G$. In fact, we prove 
that (see Lemma 2.2)
{\it if $G$ is semisimple and 
$G/L = G(p)$ is a regular orbit on $M$, with $\g$ and
$\l$  Lie algebras of $G$ and $L$ respectively, 
then the centralizer $C_\g(\l)$ is of the form 
$$C_\g(\l) = \z(\l) \oplus \a\tag1.1$$
where $\z(\l)$ is the center of $\l$ and $\a$ is either one
dimensional or it is a rank one simple Lie algebra\/}. In the following
we will distinguish these two cases by saying 
that the action of $G$ on $M$ is 
{\it ordinary\/} if $\a$ is one dimensional and
{\it extra-ordinary\/} if it
is 3-dimensional.\par
For the ordinary actions, the moment mapping $\mu$ determines a fibration
of any regular orbit of $M$ onto a flag manifold
$$\mu|_{G/L}\: G/L \to G/K\tag1.2$$
where $K = N_G^o(L)$. We say that the complex structure $J$ on $M$ is 
{\it projectable\/} if the $G$-invariant CR structure $J_{G/L}$ induced on any 
regular orbit is also $Ad(K)$-invariant and hence it descends to a
$G$-invariant  complex structure $J_K$
on $G/K$ (see \S 3). Note that 
a generic complex structure is projectable and that the 
non-projectable complex structures can occur only if the pair
$(G,L)$ verifies some special algebraic conditions (see e.g. Lemma 3.2).\par
\medskip
In this paper we obtain a classification
of the cohomogeneity one, compact K\"ahler manifolds 
with ordinary action
and projectable complex structures. We will call them 
briefly  {\it compact  K.O.P. manifolds\/} (from K\"ahler with Ordinary,
Projectable, cohomogeneity one action) and we leave the analysis
of the remaining cases to a forthcoming paper.\par
 Let  $\gamma$ be
a normal geodesic of the cohomogeneity one K\"ahler manifold
$(M, g, J)$ and for any $t\in \R$
consider the two form $\omega_t\in \Lambda^2\g$
$$\omega_t(X,Y) \= \omega_{\gamma_t}(\hat X, \hat Y)$$
where $\omega$ is the K\"ahler form of $M$ and $\hat X$,
$\hat Y$ are the infinitesimal transformations on $M$
corresponding to $X$, $Y \in \g$. 
From the semisimplicity of $G$,
it can be shown that there exists 
a unique map
$\Z\: \R \to C_\g(\l) = \z(\l)\oplus \a$ such that 
$\omega_t(X,Y) = \B([\Z(t), X], Y)$,
where $\B$ is the Cartan-Killing form of $\g$. The basic
fact of our study is that {\it 
the complex structure $J$ and the K\"ahler
metric $g$ on $M$ can be completely recovered from the complex 
structure $J_K$ induced on the flag manifold $G/K$ and
from the function $\Z(t)$\/}. \par
In order to state this result in a more precise way, let us point
out some crucial properties of $\Z(t)$ (see \S 4). 
$M$ is known to admit exactly two singular orbit 
$G/H_1$ and $G/H_2$ and we may always assume that   
$\gamma$ is parametrized so that $\gamma(0)$ and $\gamma(\d)$
are the intersections  with those  singular orbits.
Consider also an element $Z\in C_\g(\l)$
 such that 
$\a = \R Z$ and $\B(Z,Z) = -1$.
 Then, we prove that: \par
\roster
\item"a)"
the image $\ell_Z \= \Z(\R)$ is a segment parallel to $Z$
with endpoints given by $Z_1 = \Z(0)$ and $Z_2 = \Z(\d)$ and each
inner point
$P\in \ell_{Z}$ is a regular element for the flag manifold
$G/K$, i.e. $C_G(P) = K$; furthermore, if $\alpha$ is any root so 
that $E_\alpha$ is a $+i$-eigenvector of $J$, then $\alpha(P)\in i\Bbb 
R^{+}$;
\item"b)" for each endpoint $Z_i$, $i=1,2$,
$C_G(Z_i)/K \simeq \C P^{n_i}$ 
for some $n_i\geq 0$ and the complex structure $J_K$ projects naturally
onto an invariant complex structure on $G/C_G(Z_i)$.
\endroster
For a pair $(G,L)$
with $C_\g(\l) = \l + \R Z$,
we  call {\it admissible\/}
any segment which verifies  a) and b) (see Definition 5.1). 
Moreover, if $Z_i$ is one endpoint
of an admissible segment, we  call  {\it degree of
the endpoint\/} the  integer $m_i = n+1$ if  
 $C_G(Z_i)/K = \C P^{n_i}$. 
Finally, if we denote by $H_i = C_G(Z_i)$, we say that two
admissible segments $\ell_Z$ and $\ell'_Z$ are {\it
equivalent\/} if the corresponding triples
$(H_1, L, H_2)$ and $(H'_1, L, H'_2)$
determine  two equivalent abstract models
of compact cohomogenity one
$G$-manifolds with singular orbits (for the precise defintion,
see [AA]).\par
\medskip
We call {\it abstract model 
of a  K.O.P. manifold \/}
a pentuple $\Cal K = (G, L , J, \ell_Z, f)$ of the following form:
\roster
\item"i)" $G$
is a compact semisimple Lie group; $L\subset G$ is a
closed subgroup with $C_\g(\l) = \z(\l) + \R Z$; 
\item"ii)" $J$ is a 
$G$-invariant complex structure on the flag manifold
$G/K$, $K = N_G^o(L)$; 
\item"iii)" $\ell_Z$ is an admissible segment;
\item"iv)"
$f$ is a smooth real function on $\R$ with values in 
$[0,\text{length}(\ell_Z)]$, 
which is
monotone on some interval $]0,\d[$, it is invariant
by the symmetries $t\to -t,\quad t\to 2\d - t$, and
$f(0) = 0$, $f(\d) = \text{length}(\ell_Z)$,
$f''(0) = - f''(\d) = 1$.
\endroster
(see Definition 5.1).
Then our first main theorem is the following. 
\medskip
\proclaim{Theorem 1.2} Let $\Cal M$
 be the moduli space (w.r.t.
$G$-equivariant biholomorphic isometries)
of
compact {\it K.O.P.\/} manifolds  $(M, g, J)$, 
with vanishing first Betti number and let $\Cal A$ be the collection
of equivalence classes of abstract models.
Then:\par
\roster
\item there exists a bijective correspondence between
$\Cal M$ and $\Cal A$;
\item
if $(M, J, g)$ and $(M', J', g')$ are realizations
of the abstract models $\Cal K,\Cal K'$ given 
by ${\Cal K} = (G, L, J, \ell_Z, f)$
and $\Cal K' = (G', L', J', \ell'_Z, f')$ respectively, then
$M$ and $M'$ are equivalent as $G$-manifolds if and 
only if $G/L = G'/L'$ and $\ell_Z \simeq \ell_Z'$.
\endroster
\endproclaim
For the exact relation between an abstract model and the geometric 
data of its realizations, see Cor. 4.5 and Thm. 5.2.
\bigskip
Let us come now to the  K\"ahler-Einstein manifolds. From Theorem 1.2,
the classification of compact K\"ahler-Einstein manifolds
reduces to  characterizing the abstract models associated 
to Einstein metrics. The result is easily obtained 
and it is the following second main theorem. \par
In the statement, we denote by 
$\g = \l + \R Z + \m$ the orthogonal decomposition w.r.t. $\B$
of the Lie algebra $\g$ of $G$ with   $\l + \R Z = C_\g(\l)$.
By $R$ we denote the roots of $\g^\C$ corresponding to a 
Cartan subalgebra $\h \subset \k^\C$  and  
by $R_\m$ the roots corresponding to root vectors in $\m^\C$. We 
also assume that $R^+_\m$ stands for the  roots so that
$J_K(E_\alpha) = + i E_\alpha$ and
$Z^\kappa$ denotes the element
$Z^\kappa = \sum_{\alpha \in R^+_\m} H_\alpha$, where 
$H_\alpha$ denotes the element of $\h$ corresponding to the 
root $\alpha$ under the isomorphism $\B^*:\h^*\to \h$ induced 
by the Cartan-Killing form $\B$. 
\par
\bigskip
\proclaim{Theorem 1.3} An abstract model  $\Cal K =(G, L, J, \ell_Z,
f)$ corresponds to an Einstein-K\"ahler manifold with Einstein 
constant $c=1$ if and
only if the following conditions holds:
\roster
\item if $Z_1$ and $Z_2$ are the two endpoints of $\ell_Z$ and $m_1$
and $m_2$ are their degrees,  then 
$$|Z_1 - Z_2| = m_1+m_2\qquad,\qquad
\frac{m_2}{m_1+m_2} Z_1 + \frac{m_1}{m_1+m_2} Z_2 = Z^\kappa\ ;\tag1.3$$
\item the following Futaki integral vanishes
$$\int_{-m_1}^{m_2}y\prod_{\alpha\in R^+_\m}
\alpha(Z^\kappa - y Z)dy = 0\ ; \tag1.4$$
\item $f(\d) = m_1+m_{2}$ and the inverse function $t(f)$ of $f|_{[0,\d]}$
is 
$$t(f) = \int_{0}^f \frac{\sqrt{\prod_{\alpha\in R^+_\m}
\alpha(Z^\kappa- s)Z}}
{
 - 2
 \int_{-m_1}^s v {\prod_{\alpha\in R^+_\m} \alpha(Z^\kappa- vZ)}
 dv
}
ds\ .\tag1.5$$
\endroster
A realization $(M,J)$ of an abstract model $\Cal K$
has two singular orbits of codimension $2m_1$ and $2m_2$,
respectively, and it
admits a K\"ahler-Einstein metric with Einstein constant $c=1$ if and only if
the admissible segment $\ell_Z$ verifies (1) and (2); on such 
complex manifold the K\"ahler-Einstein
 metric
is unique.  
\endproclaim
\bigskip
\remark{Remark 1.4} All the examples 
of compact K\"ahler-Einstein metrics
given by Koiso and Sakane in [KS] have ordinary 
action and  projectable complex structure and 
they exhaust the realizations of the abstract models
of Theorem 1.3.
Therefore, our result can be
restated saying that {\it  the
 compact, K\"ahler-Einstein K.O.P. manifolds
are exactly all those given by Koiso and Sakane in [KS]\/}.
\footnote {We like to mention that, after a preliminary version of this paper
([PS]) 
had been completed, we received an interesting preprint by A. Dancer and M. Wang
([DW]) concerning 
non compact K\"ahler-Einstein metrics of cohomogeneity one. }
\endremark
\bigskip
\remark{Remark 1.5} The admissible segments with
$m_1 = m_2 = 1$, which verify the condition
(1) of Theorem 1.3 are diameters of the sphere
$\Cal S = \{ Z\: |Z - Z^\kappa| = 1\}$. If $\Cal S$
is included in the $T$-chamber of $Z^\kappa$, it is not hard
to select at least one admissible segment which verifies
also condition (2).  This gives a way to build 
compact K.O.P. manifolds, with singular
orbits of codimension 2 and   which admit a K\"ahler-Einstein metric.\par
The admissible segments with $m_1\neq 1$ or $m_2\neq 1$, which verify
(1)  of Theorem 1.3 
 are much more rare, so that 
condition (2) of the same
theorem creates  an extremely strong obstraction
to the existence of  K.O.P. manifolds with Einstein metric
and singular orbits with codimension higher than $2$. Examples
of this kind are given by the complex projective
space $\Bbb {CP}^n$: the group $G\subset U(n+1)$ can be taken either to have 
one fixed point or to be a product group $G=G_1\times G_2$, 
where the semisimple parts $G_i^{ss}=SU(n_i)$
or $Sp({{n_i}\over {2}})$ for $i=1,2$ and $n_1+n_2=n+1$. 
It would be interesting to know other non trivial 
examples of such K.O.P. manifolds and if they can be 
explicitly listed.
\endremark
\bigskip
We conclude observing that any compact K.O.P. manifold is also 
almost homogeneous w.r.t. to the action of the complexified group 
$G^{\C}$; we refer to [HS] for a general theory.
\vfil\eject
\subhead 2. First properties of compact cohomogeneity
one K\"ahler manifolds \endsubhead
\bigskip
\subsubhead 2.1 Cohomogeneity one compact K\"ahler
manifolds
\endsubsubhead
Throughout the following $(M,g,J)$ will denote a compact K\"ahler
manifold with vanishing first
Betti number, where 
$g$ is the metric tensor, $J$ the complex structure and
$\omega = 
g(\cdot,J\cdot)$ is the associated K\"ahler fundamental form.
The  Lie algebra of a Lie
group acting on $M$
will be always denoted  by the corresponding gothic letter.
Furthermore, if $X$ is any element of $\frak g$, 
the symbol $\hat X$ will be used for the corresponding
 Killing vector field on $M$. Finally,
$\B$ denotes the Cartan-Killing form of any
semisimple Lie algebra.\par
We assume that $G$ is a compact  connected
 Lie group of isometries of $M$ (and hence also 
of holomorphic transformations, since $M$ is compact) with at least one
hypersurface orbit. In this situation, it is well known (see [Br],
for instance) that either all orbits are regular and of codimension
one or they are all regular except exactly two singular orbits. The
orbit space $\Omega = M/G$ is  diffeomorphic to  $S^1$ 
in the first
case, while it is equivalent to $[0,2\pi]$ in the second case. 
Indeed, assuming that $b_1(M) = 0$, there is no fibration
$M\to S^1$ and so, {\it $M$
admits exactly two singular orbits $S_1$ and $S_2$\/}.
The   subset $M_{\text{reg}}$ of all regular points 
is dense in $M$ and the stability subgroups 
are all conjugate to the same compact subgroup, say
 $L$.\par
Since the manifold $M$ is orientable, every regular 
 orbit is orientable and we may define a unit normal
vector field $\xi$ on the whole subset $M_{\text{reg}}$. 
It is known (see e.g. [AA1]) that any
integral curve of $\xi$ in $M_{\text{reg}}$ is a geodesic. Actually,
for any point $p \in M_{\text{reg}}$, the geodesic 
$\gamma:\Bbb R\to M$ such that 
 $\gamma'_0 = \xi(p)$  intersects orthogonally 
all regular orbits and it is called {\it normal geodesic\/}.
It is clear tha for any $t$, $\gamma'_t$ is
equal or opposite to $\xi_{\gamma_t}$. More precisely,
assume that $\gamma$ is parametrized so that 
$A=\gamma^{-1}(M_{\text{reg}}) \subset \R$ is of the form
$A= \bigcup_{k\in\Bbb Z}]k\d,(k+1)\d[$,
where $\d={\text{dist}}(S_1,S_2)$.
Then,  up to reversing the 
orientation, it is easily seen that 
$$\gamma'_t = (-1)^k \xi_{\gamma_t}\quad t\in ]k\d,(k+1)\d[.\tag2.1$$
\bigskip
\subsubhead 2.2 The moment mapping and the flag manifolds associated
to the orbits
\endsubsubhead We briefly recall the definition of the {\it moment mapping\/}.
For a complete discussion of its main properties
see e.g. [GS].\par
For any $X\in \g = Lie(G)$ the facts  that
$\Cal L_{\hat X}\omega = d(\imath_{\hat
X}\omega)=0$ and that $b_1(M)=0$ imply the existence a potential 
$f_X\in C^\infty(M)$ such that $df_X=\imath_{\hat X}\omega$. 
By the compactness of $M$, it can be proved   the potential 
$f_X$  can be chosen for all $X$ in a consistent way
so that the
mapping $\mu:M\to \g^*$ 
$$\mu(m)(X) = f_X(m),\quad m\in M$$
is well defined and $G$-equivariant (we consider the coadjoint
action of $G$ on $\g^*$). Any such map  is 
usually called {\it moment mapping\/} and it is unique 
in case $G$ is semisimple.
\par
The existence of a moment mapping $\mu$ brings to the following
basic result.\par
\medskip
\proclaim{Proposition 2.1} Suppose that $G$ is a compact connected
Lie group acting on  $M$ with cohomogenity
one. 
Let $\gamma\: \R \to M$
be a normal geodesic, $L\subset G$ be the common stability
subgroup at all the points $\gamma(\R) \cap M_{\text{reg}}$ and,
for any $t$, let 
$$K_t = \{g\in G;\ Ad(g)^*(\mu(\gamma_t))=\mu(\gamma_t)\} 
\supset G_{\gamma_t}\ .$$
Then
\roster
\item $G/K_t$ is a flag manifold for any $t\in \R$;
\item  for $t\in A =\gamma^{-1}(M_{\text{reg}})$, 
$L$ is a subgroup 
of codimension one in $K_t$. 
\endroster
Furthermore, if $G$ acts almost effectively, then the center 
$Z(G)$ is at most one dimensional.
\endproclaim
\demo{Proof} 
(1) follows from the fact that $\mu(G\cdot \gamma_t) =
G/K_t$  is a coadjoint orbit in $\g^*$ (for this
property and others regarding the flag manifolds,
see [Be]). For (2)
we only need to evaluate
the dimension of the fiber of the 
mapping $\mu_p\: G/L \to G/K_t$ for any $p=\gamma_t$ with $t\in A$.\par
If $v\in \ker d\mu_p\cap T_p(G/L)$, we may find $Y\in
\g$ with $\hat Y_p=v$ and we have, for all $Z\in {\frak g}$
$$\eqalign {0 =&\left.{{d}\over{ds}}\right|_{s=0}\mu(\exp(s{Y})p)(Z) =
\left.{{d}\over{ds}}\right|_{s=0}f_Z(\exp(s {Y})p) =\cr 
{} &= df_Z|_p(Y) = (i_{\hat Z}\omega)|_p(Y) = \omega_p(Z,Y).\cr}$$
Since the pull back of $\omega$ to the orbit $G\cdot p = G/L$ 
has one-dimensional
kernel, it follows that  $\dim\ker d\mu_p = 1$, and hence 
$\dim K_t/L=1$.\par
Finally, we note that if $G$ acts (almost)
effectively on $M$, then so does it on each regular orbit $G/L$; if
$\dim Z(G)>1$, since $Z(G)\subset K_t$ for all $t$, we would have that
$\dim Z(G)\cap L >0$ and this contradicts the (almost) effectiveness of the
$G$-action on $G/L$.
\qed
\enddemo 
\bigskip
From the previous result, 
 $G$ 
is either semisimple or it admits a one dimensional center.
We will discuss the second case in a separate paper
and from now on, we will suppose that 
{\it the group $G$ is semisimple\/}. Recall that
being $G/K_t$ a flag manifold, $K_t$
is a maximal rank subgroup of $G$. 
This  puts strong restrictions on the pair 
$(G,L)$. In fact
\medskip
\proclaim{Lemma 2.2} Let $G$ be a compact, semisimple Lie group 
$G$ and let $L$ be a closed Lie subgroup which is a codimension one 
subgroup of a maximal rank Lie subgroup $K$ of $G$. 
Then the centralizer
$C_\g(\l)$ of the Lie algebra $\l$ in $\g$ admits a decomposition 
$$C_\g(\l) = \z(\l) \oplus \a$$ 
where $\z(\l)$ is the center of $\l$ and $\a$ is either a 1-dimensional
Lie subalgebra or a rank one simple Lie algebra with
$\a^\C = \g(\alpha) = span_\C<E_{\pm\alpha}, H_\alpha>$, for some root
$\alpha$ of a Cartan subalgebra  $\t^\C \subset\k^\C$.
 \endproclaim
\demo{Proof} Let $\l$ and $\k$ the Lie algebras of $L$ and $K$
and let $\l^\perp$ the orthogonal complement of $\l$ in $\g$
w.r.t. the Cartan-Killing form $\B$. By hypotheses,
we may decompose the Lie algebra $\k$ of $K$ as 
$$\k = \l +\R Z,$$
for some non zero $Z\in(\k\cap \l^\perp)$.
Obviously  $[\l,Z]=0$. Moreover, if we fix a
maximal torus $\t$ of $\l$, we have that $\hat \t=\t+\R Z$ is a
maximal torus for $\k$ and hence also
for $\g$. It then follows that 
$\hat \t^\C$ is a Cartan subalgebra of
$\g^\C$ and we have the standard decomposition  
$$\g^\C = \hat \t^\C +
\sum_{\alpha\in R}\C E_\alpha\ ,$$ 
where $R$ denotes the root system of
the pair $(\g^\C,\hat \t^\C)$. Moreover we have 
 $$\k^\C = \l^\C + \C Z = \t^\C +
\sum_{\beta\in R_K}\C E_\beta,$$ for some root subsystem $R_K\subset R$. \par
Now, it is clear that the centralizer 
$C_{\g}(\l)$ contains $\z(\l)+\R Z$ and, if 
it is strictly bigger, then there is at least one 
$\alpha\in R\setminus R_K$
with  $E_\alpha\in C_{\g^\C}(\l^\C)$. 
It then follows that $\alpha$ vanishes on $\t^\C$; 
since $\t^\C$ has codimension one in the CSA, 
there are at most two roots, say 
$\{\pm \alpha\}$, with this property and the centralizer 
$C_{\g^\C}(\l^\C)$ is given by
$$C_{\g^\C}(\l^\C) = \z(\l^\C) + \a^\C,\quad \a^\C =  
span_{\C}<Z,E_\alpha,E_{-\alpha}>.$$ 
Note that $[E_\alpha,E_{-\alpha}]$ is $\Cal B$-orthogonal to $\t^\C$, 
hence is a multiple
of $Z$ and we have that $\dim C_G(L) - \dim Z(L) =3$. In particular,
the 
sum $\z(\l)+ \a$ coincides with $C_\g(\l)$ and it
is a direct sum.
\qed
\enddemo
\bigskip
 In the hypotheses of Lemma 2.2, it can be proved also the 
follwing fact (see [AS]): {\it
 $\a$ is a rank one
simple Lie algebra if and only if $G$ is simple, 
$L = C_G(\g(\alpha))$ and the root  $\alpha$ verifies the following
property:\ \rm for any root $\beta$ which is orthogonal to 
$\alpha$, $\alpha\pm \beta \notin R$\it\/}. \par
In particular, 
from this it can be checked
that  the pairs $(G,L)$, with $G\neq G_2$ and 
for which $\a$ is a 3-dimensional simple Lie algebra
are in one to one 
correspondence with the 
homogeneous quaternionic K\"ahler manifolds of simple
Lie groups (the so called {\it Wolf spaces\/}).\par
\bigskip
\definition {Definition 2.3} 
Let $G$ be a compact semisimple Lie group and let $L$ be 
a compact subgroup which is contained 
in a subgroup of maximal rank as a subgroup of 
codimension one. We will say that the pair $(G,L)$ is {\it 
ordinary \/} if 
$$\dim C_G(L) = \dim Z(L) + 1$$
where $C_G(L)$ and $Z(L)$ are the centralizer of $L$ in $G$ 
and the center of $L$, respectively. It is called 
{\it extra-ordinary\/}
otherwise.
Similarly, we will also 
say that the action of  $G$ on $M$ is {\it ordinary}
(resp. {\it extra-ordinary})  
if the regular orbits are of codimension one and 
$G$-equivalent to $G/L$, for some 
ordinary (resp. extra-ordinary) pair $(G,L)$.
\enddefinition
\bigskip
 We conclude this section with the 
following immediate Corollary of Proposition 2.1.
\bigskip
\proclaim {Corollary 2.4} 
Let $G$ be a compact semisimple Lie group acting on $M$ 
with ordinary  action and let $\gamma$ be a normal
geodesic. Then for any $t\in A = \gamma^{-1}(M_{\text{reg}})$, 
the moment mapping $\mu$
fibers the orbit
$G\cdot \gamma_t = G/L$ 
over the same flag manifold $G/K$.
\endproclaim
\bigskip
\bigskip
\subhead 3. The invariants of cohomogeneity one compact 
K\"ahler manifolds
\endsubhead
\bigskip
The purpose of this section is to determine a set 
of invariants which can be defined on the family
of flag manifold $G/K_t$ of Proposition 2.1 and by which
the   complex  and the Riemannian structure
of $M$ can be uniquely characterized. \par
We keep the same notations as in the previous section;
we also fix the following
$\Cal B$-orthogonal decomposition of  $\g$ 
$$\g = \l+\n = \l + \a + \m\ ,\tag3.1$$ 
where $\n = \a + \m$ and $\a$
is the orthogonal complement to $\z(\l)$ in $C_\g(\l)$, as in Lemma
2.2. In case the action of $G$ is ordinary, we set
$\a = \R Z$ for a fixed element $Z$ with ${\Cal B}(Z,Z) = -1$.\par
By means of the decomposition
(3.1),  for each $t\in A = \gamma^{-1}(M_{\text{reg}})$, we
have an identification map  
$\phi_t:\n\to T_{\gamma_t}(G\cdot\gamma_t)$ given by 
$$\phi_t(X) = \hat X_{\gamma_t}\ .\tag3.2$$
\bigskip
\subsubhead 3.1 The invariant associated to the complex structure
\endsubsubhead
The complex structure $J$ of $M$ 
induces an integrable $G$-invariant CR-structure on each 
regular orbit $G/L$. For any $p\in G/L$,  we denote by
$\D_p$  the maximal $J_p$-invariant 
subspace of $T_p(G/L)$ so that the CR structure on $G/L$
 is uniquely determined by the pair $(\D, J|_\D)$.\par
By means of the identification (3.2), for any
$t\in A$  we  define the
ad($\l$)-invariant 
subspace $\m_t\subset\n$ 
$$\m_t = \phi_t^{-1}({\Cal D}_{\gamma_t})\ .\tag3.3$$
Note that for each $\m_t$ there exists a unique (up to sign)  
ad($\l$)-invariant 
element 
$Z_\D(t) \in \n$ so that
$$Z_\D(t) \in \n \cap \m_t^\perp\ ,\qquad \B(Z_\D(t), Z_\D(t)) = -1\ .\tag3.4$$
Clearly, $Z_\D(t) \in C_\g(\l)$ and hence, if the action of $G$ is ordinary,
 we may assume that
$$Z_\D(t) \equiv Z\ .\tag3.5$$
By (3.2) we also define a complex 
structure $J_t$ on $\m_t$ for any $t\in A$
 $$J_t=\phi_t^*(J|_{{\Cal D}_{\gamma_t}})\ .$$
In the following, we will  consider $J_t$ as trivially 
extended to an endomorphism of $\n$ for any $\t$.
It is clear that each $J_t$ is ad($\l$)-invariant. 
From the integrability of the CR-structures $(\D, J|_\D)$
it can also easily proved that, for all 
$X,Y\in\m_t$, 
$$[J_tX,Y]_{\n}+[X,J_tY]_{\n}\in\m_t\tag 3.6$$
$$J_t([J_tX,Y]_{\n}+[X,J_tY]_{\n}) = ([J_tX,J_tY]-[X,Y])_{\n}.\tag 3.7$$
Finally, for any $t\in A$,
we may consider  the unique (up to sign) element  
$$Z_{J\xi}(t) \in \R \phi_t^{-1}(J\xi_{\gamma_t})\ ,
\qquad
\B(Z_{J\xi}(t),Z_{J\xi}(t)) =-1 \ .\tag3.8$$
We note that even the vector $Z_{J\xi}(t)$ is in $\a$, because the
tangent vector $J\xi_{\gamma_t}$ is left fixed under 
the isotropy representation of $L$. Therefore 
 if the action of $G$ is ordinary, we may always choose
$$Z_{J\xi}(t) \equiv Z\ .\tag3.9$$
It is immediate to realize that  the map
$${\Cal J}\: A \to \a\times End(\n) \times \a$$
$${\Cal J}(t) \= (Z_\D(t), J_t, Z_{J\xi}(t))$$
determines uniquely the complex structure $J$ at all tangent spaces
$T_{\gamma_t}M$ and, being $J$ $G$-invariant, it characterizes
completely $J$ over $M_{\text{reg}}$. Furthermore, it
is unique up to $G$-equivalence (i.e. $G$-invariant biholomorphic isometry)
 in the following sense: {\it if $(M,J,g)$ and
$(M',J',g')$ are $G$-equivalent and $\gamma$ and
$\gamma'$ are two normal geodesics of $M$ and $M'$
which are mapped one onto the other by a $G$-equivariant
isometric biholomorphism, then the corresponding maps
${\Cal J}$ and ${\Cal J}'$ coincide (up to the signs of $Z_\D$ and
of $Z_{J\xi}$)\/}.\par
We will call ${\Cal J}$ the {\it invariant associated to 
the complex structure\/}.\par 
\medskip
Finally, let us consider the following definition: we say that the complex
structure $J$ of $M$ is {\it projectable\/} if 
$$Ad(K_t) \cdot J_t = J_t\tag3.10$$
for any $t\in A$. This  means that each CR structures
$J_t$ descends to
a $G$-invariant almost complex structure $J_{K_t}$
on the flag manifold $G/K_t$
 and such that the natural map
$\pi\: G/L \to G/K_t$ verifies
$$\pi_*(J_t(X))  = J_{K_t} \pi_*(X)\qquad \forall X\in \D$$
Note that  (3.7) implies that each $J_{K_t}$ 
 is an {\it integrable\/}
complex structure on $G/K_t$. \par
\medskip
Now let us limit ourselves to considering ordinary actions. In this
case, we have already seen that 
$Z_\D(t) \equiv  Z_{J\xi}(t) \equiv Z$ and hence 
that ${\Cal J}$ is uniquely given by the family of complex structures
$J_t$ on the subspace $\m\equiv\m_t$. If furthermore
$J$ is projectable,   each $J_t$ is mapped onto 
an invariant  complex structure of the same flag manifold
$G/K_t \equiv G/K$ (by Corollary 2.4).
Since the invariant complex structures
on a flag manifold are a discrete finite set (see e.g. [Be]), we
 have that the complex
structures $J_t$ are constant on each connected
interval of $A$. Hence we conclude that
{\it if the action of $G$
is ordinary and $J$ is projectable, then the invariant $\Cal J$
is uniquely determined by a corresponding  
$G$-invariant complex structure $J_K$ on the 
flag manifold $G/K$\/}.\par
\medskip
\remark{Remark 3.1} It can be proved that,
for a {\it generic\/} compact K\"ahler manifold
with ordinary, cohomogeneity one action, the complex
structure $J$ {\it is\/} projectable (see [AS]). However, a useful
criterion which guarantees the projectability of the 
complex structure is given by the following 
lemma of straightforward proof.
\endremark
\bigskip
\proclaim{Lemma 3.2} Suppose that the action of $G$ on $M$
is ordinary and consider the decomposition 
$\g = \l +  \R Z + \m$ in (3.1). If every
ad($\l$)-irreducible submodule of $\m$ appears with multiplicity one,
then the complex structure $J$ is projectable.
\endproclaim
\bigskip
\example{Example 3.3}  Consider the Hermitian symmetric 
space $M=SO(n+2)/SO(2)\times SO(n)$, which 
is the complex quadric $Q_n$ of ${\Bbb CP}^{n+1}$ 
and which provides one compactification of the 
tangent bundle of a sphere $S^{n+1}$. 
The group $G=SO(n+1)$ acts with cohomogeneity one on $M$ with 
principal isotropy subgroup $L=SO(n-1)$. A normal geodesic 
is explicity given by 
$$\gamma_t = \left[\left(\matrix\cos t& 0& -\sin t & \cr
0&1&0&\Bbb O\cr \sin t&0&\cos t& \cr &\Bbb O& & \Bbb I_{n-1}
\endmatrix \right)\right],$$
where $[\cdot]$ denotes the projection of an element of $SO(n+2)$ 
onto the quotient space $M$.
The pair ($G,L$) is ordinary, but 
the complex structure is not projectable. In fact,
  keeping the notation as above, 
if we write $\g =\l+\Bbb R Z+\m$, then $\m$ splits as the 
sum of two equivalent, irreducible 
ad($\l$)-submodules $\m_1,\m_2\cong \Bbb R^{n-1}$; 
moreover an easy computation shows that the complex
structure $J_t\in End(\m)$ is given, for $(v_1,v_2)\in \m_1+\m_2$, by
$$J_t((v_1,v_2)) = (v_2/\sin t,-(\sin t) v_1).$$
These endomorphisms do not commute with ad($Z$) and hence with
ad($\k$). Note that ad($Z$) represents 
the standard ad($\l$)-invariant
complex structure $(v_1,v_2)\mapsto (-v_2,v_1)$ on $\m$.\endexample
\bigskip
\subsubhead 3.2 The invariants associated to the K\"ahler form
and to the metric tensor
\endsubsubhead
Let us consider the usual normal geodesic $\gamma$. For any $t\in \R$,
consider the following 2-form $\{\omega_t\}_{t\in\R}$ 
on  $\g$: 
$$\omega_t(X,Y) \= \omega_{\gamma_t}(\hat X,\hat Y), \quad X,Y\in\g.$$
Since the K\"ahler form $\omega$
is closed and every Killing vector field preserves it, we have that
$$0=3d\omega(\hat X,\hat Y,\hat W) = \omega(\hat X,[\hat Y,\hat W]) + 
\omega(\hat W,[\hat X,\hat Y])
+ \omega(\hat Y,[\hat W,\hat X])\ ,$$
so that for all $t\in\R$ and for all $X,Y,W\in\g$ we have 
$$\omega_t([X,Y],W) + \omega_t([W,X],Y) + \omega_t([Y,W],X) = 0\ .
\tag 3.11$$
Since $\g$ is semisimple, we may find $F_t\in End(\g)$ with 
$\omega_t(X,Y) = \B(F_t(X),Y)$ for all $X,Y\in\g$. Now, from
(3.11)  and  the non degenerancy of $\B$, we get  that
$$F_t([X,Y]) = [F_t(X),Y] + [X,F_t(Y)],$$
that is each $F_t$ is a derivation of $\g$. But then
 for each $t \in\R$ there exists an 
unique $\Z$ such that 
$$F_t = ad(\Z(t))\ .$$
Note that,
for all $t\in A$, we may say that 
$\ker F_t = \l + \R v_t$, 
where $v_t\in \l^\perp$ and 
$\hat v_t|_{\gamma_t} = J\xi$. 
Therefore $\ker F_t$ coincides with the Lie algebra 
of the subgroup $K_t$ defined in Proposition 2.1
and  the flag manifold $G/K_t$ is $G$-equivalent 
to the $G$-orbit of  $\Z(t)$ in $\g$. 
In particular, for any $t\in A$,
$\Z(t)$ is $ad(\l)$-invariant and $\Z(t) \in \z(\l) + \a$.\par
In conclusion,
to the $G$-equivalence class of $(M,g,J)$ we may uniquely associate
the map
$$\Z\: \R \to \g\qquad, \qquad  \Z|_A\: A \to \z(\l) + \a$$
$$\B([\Z(t),X],Y)  = \omega_{\gamma_t}(\hat X,\hat Y)\tag3.12$$
We will call $\Z(t)$ {\it the invariant associated to the
K\"ahler form\/}.\par
\bigskip
\remark{Remark 3.4}
If we decompose $\g=\k_t+\m_t$ with $\m_t$ $\B$-orthogonal and 
ad($\k_t$)-invariant subspace,  the restriction 
$\omega_t|_{\m_t}$ is
ad($\k_t$)-invariant and hence it defines a symplectic form on 
the flag manifold $G/K_t$. \par
Note that, by the previous observations,
$\Z(t)\not=0$ for all $t\in A$; it can be checked
that $\Z(t)=0$ only when the orbit
$G\cdot\gamma_t$ is a point or a totally real submanifold.
\endremark
\bigskip
Now, let us use the identification (3.2) to determine an invariant 
which can be used to characterize
 the metric tensor. For each $t\in A$, let us consider the
ad($\l$)-invariant, positively definite bilinear form 
on $\n$ given by the pull back of the Riemannian 
metric:
$$g_t\=\phi_t^*(g|_{T_{\gamma_t}G(\gamma_t)}).$$
Note that $g_t|_{\m_t}$
is completely determined 
by the invariants $\Cal J$ 
and $\Z$: in fact for all $t\in A$
$$g_t|_{\m_t} = -\phi_t^*(\omega(\cdot, J \cdot)|_{\D_{\gamma_t}}) =
-\B([\Z(t), \cdot], J_t\cdot)|_{\m_t}$$
 Therefore, once $\Cal J$ and $\Z$ are given, the only  ingredient
that is necessary to recover the whole tensor $g_t$
is the function 
 $g_t(Z,Z) = ||\hat Z_\D(t)||^2$.
We define as {\it invariant
associated to the metric tensor\/}   the following map 
$a\: \R \to \R$
$$a(t) = \left\{\matrix (-1)^{k}||\hat Z_\D(t)||_{\gamma(t)} & t 
\in ]k\d, (k+1)\d[\\
0 & t \in {\Bbb Z} \d
\endmatrix \right.\tag3.13$$ 
where $\d$ is the distance
between the two singular orbits and
the parametrization of $\gamma$ is assumed to satisfy (2.1).
The definition has been motivated by the fact that,
if the action is ordinary and
$Z_\D \equiv Z$, such a function is indeed ${\Cal C}^\infty$
over $\R$, as we will show in the next section.\par
\bigskip
\bigskip
\subhead 4. The properties of the invariants $\Cal J$, $\Z$ and $a$
on the K.O.P. manifolds
\endsubhead
\bigskip
From now on, we will assume that the action of $G$ is ordinary
and that the complex structure is projectable, i.e. that $M$
is a compact {\it K.O.P.\/} manifold. The first
fact is that $\Cal J \equiv (Z, J, Z)$
where $J$ is a CR structure which is projectable onto an
invariant complex structure on $G/K$. We continue with the following
\bigskip
\proclaim{Proposition 4.1} Let $\Z$ and $a$ be the invariants associated
to the K\"ahler form and the metric tensor of $M$. Then
they  are both smooth  functions of $t$ and they
verify the ordinary differential equation
$${d\over{dt}}\Z(t)= - a(t)Z\ .
\tag4.1$$
Moreover, if we denote by $H_1,H_2$ the two singular stability subgroups 
of the singular points ${\gamma_0}$ and ${\gamma_{\d}}$ respectively,
then:\roster
\item $H_1\cap H_2 = K$; 
\item 
the singular orbits $G/H_i$ are complex submanifolds of $M$; 
\item $H_i = C_G(\Z(t_i))$, $t_i = 0$ or $\d$, and 
$\mu(G/H_i)\simeq G/H_i$;
\item  any homogeneous space 
$H_i/K$ is diffeomorphic to a complex projective space;
\endroster 
In particular the geodesic symmetries $\sigma_1,\sigma_2$ at 
$\gamma_0$ and $\gamma_{\d}$ belong to $K$. 
\endproclaim
\demo{Proof}
For  the proof, we only need to show that (4.1) is verified at all 
$t\in \R$, so that the smoothness of $a$  will
follow immediately from the smoothness of $\Z$.\par
In order to do this, we
 write down the K\"ahler condition $d\omega=0$ and find the corresponding
condition on $\Z$.
Let us consider two arbitrary vectors $v_1,v_2 \in\m$, so 
that $\phi_t(v_i)\in{\Cal D}_{\gamma_t}$ for all $t\in A$. 
The condition $d\omega=0$ imply that 
$$\eqalign { 0&= 3d\omega(\xi,\hat v_1,\hat v_2) = 
\xi\omega(\hat v_1,\hat v_2) - \hat v_1\omega(\xi,
\hat v_2) + \hat v_2\omega(\xi,\hat v_1) \cr
{}&-\omega([\xi,\hat v_1],\hat v_2) + \omega([\xi,\hat v_2],\hat v_1) - 
\omega([\hat v_1,\hat v_2],
\xi)\cr
{}& = \xi\omega(\hat v_1,\hat v_2) - \omega(\xi,[\hat v_1,\hat v_2]) + 
\omega(\xi,[\hat v_2,
\hat v_1]) - \omega([\hat v_1,\hat v_2],\xi)\cr
{}& = \xi\omega(\hat v_1,\hat v_2) + 
\omega([\hat v_1,\hat v_2],\xi),\cr}\tag 4.2$$ 
where we have
used the fact that ${\Cal L}_{\hat v_i}\omega=0$ and 
$[\hat v_i,\xi]=0$ for $i=1,2$.
Now (4.2) reduces to the following for $t\in ]k\d,(k+1)\d[$
$$(-1)^k{d\over{dt}} \B([\Z(t),v_1],v_2) - g_t([v_1,v_2]_\n,
(-1)^k\frac{1}{a(t)} Z) = 0.\tag 4.3$$
Since now $Z$ and $\m$ are both $g_t$- and $\B$-orthogonal, 
(4.3) becomes
$$(-1)^k{d\over{dt}} \B(\Z(t),[v_1,v_2]) + (-1)^k\frac{1}{a(t)}
g_t(Z,Z)\B([v_1,v_2]_\n,Z) = 0,$$
or equivalently 
$$(-1)^k{d\over{dt}} \B(\Z(t),[v_1,v_2]) + (-1)^k {a(t)}\B([v_1,v_2],Z) = 0,$$
so that the vector ${d\over{dt}}\Z(t)+ {a(t)}Z\in\z(\k)$ 
must be $\B$-orthogonal to 
$[\m,\m]$. We now recall the decompositions
$$\k^\C = \t^\C + \sum_{\alpha\in R_K}\C E_\alpha,\quad 
\g^\C = \t^\C + \sum_{\alpha\in R}\C E_\alpha = \k^\C + 
\sum_{\alpha\in R\setminus R_K}\C E_\alpha,$$
where $\t$ is the Lie algebra of a maximal torus of $K$, 
$R$ is the root system of the pair 
$(\g^\C,\t^\C)$ and $R_K$ is the root subsystem given by 
$R_K = \{\alpha\in R;\alpha|_{\z(\k)} =
0\}$. Now if $v={d\over{dt}}\Z(t)+a(t)Z\not=0$, we choose a root 
$\beta\in R$ such that
$\beta(v)\not=0$; then $\beta\not\in R_K$ and therefore
$E_\beta,E_{-\beta}\in\m^\C$; so  
$$0\not= \beta(v) = \B(v,[E_\beta,E_{-\beta}])\in \B(v,[\m^\C,\m^\C]).$$
This contradiction shows that $v=0$ and equation (4.2) is equivalent to (4.1)
at all $t\in A$.\par
It remains to show that (4.1) holds also  when $t\in {\Bbb Z} \d$.
 By definition of $a$, this amounts to show that 
$a$ is everywhere continuous, or that $\hat Z$ vanishes at the singular points
$\gamma_{k\d}$ for $k\in\Bbb Z$.
We show this for $t=0$, since for all other $t = k \d$, the argument
would be the same. \par
Let $H_1\supset L$ be the stability sugroup of the singular point 
$\gamma_0$ and let $\h_1$ be its Lie algebra. Consider the decomposition
$\h_1 = \l + \p$, with $\p=\h\cap\l^\perp$ and
$[\l,\p]\subseteq\p$. \par
Assume that  $Z\not\in\h_1$.
Then
$Z\not\in\p$, so that $Z=Z_1+Z_2$ with $Z_1\in\p$ and 
$0\not=Z_2\in\p^\perp\cap\l^\perp$. Since $Z$
centralizes $\l$, we have that $Z_2\in C_{\g}(\l)\cap\l^\perp$
and hence $Z_2\in\Bbb R Z$. Therefore we conclude that $Z\in
\p^\perp$ and that $\p\subset\m$. \par
We now take $v\in\p\setminus\{0\}$ and consider the
function $\varphi(t)=||\hat v||^2_{\gamma_t}$ with $t\in \Bbb R$. 
Since $v\not\in\l$, the function
$\varphi$ has a minimum for $t=0$, so that  $\varphi'(0)=0$.
 On the other hand, for $t\in ]0,\d[$,
we have that  
$$\varphi(t)=\omega(J\hat v,\hat v)_{\gamma_t} = \B([\Z(t),Jv],v)\ ,$$
and therefore, by (4.1) on $]0,\d[$,
$$\der \varphi(t) = -||\hat Z||_{\gamma_t} \B([Z,Jv],v)\ ,$$
with $\B([Z,Jv],v)\not=0$, since the function $\varphi$ is not constant. 
This implies that 
$\hat Z_{\gamma_0}=0$, hence $Z\in\h_1$. This contradiction concludes 
the proof of (4.1). Furthermore we also obtained that 
$K\subseteq H_1\cap H_2$.\par
Next, we prove (2), i.e. the singular orbits are complex submanifolds.\par
We recall the decomposition 
$\g = \k + \m$ and observe that, being $G/K$ a flag
manifold,  the $\k$-module $\m$ splits as 
sum of irreducible and inequivalent $\k$-submodules $\m_i,i=1,\dots,s$,
each of them  stable under the projectable complex structure $J$.
Since $\k\subseteq\h_i$, then the $\B$-orthogonal
complement $\h^\perp_i$ in $\g$ is sum of some of 
the submodules $\m_j$ and therefore is
$J$-stable. Let $v\in\h_i^\perp$ and consider the vector field 
$\hat v_{\gamma_t}$ for $t\in
[0,\d[$; along the normal geodesic $\gamma_t$ with $t\in ]0,\d[$, we have 
$$J\hat{v}_{\gamma_t}=\widehat{Jv}_{\gamma_t}.$$
Since both members are continuous vector fields, we get that 
$J\hat v_p = \widehat{Jv}_p\in T_pG/H_i$,
where $p=\gamma_0$ and this proves that each singular orbit is complex.\par
It is now easy to recognize that each $\goth h_i$ is the 
centralizer in $\g$ di $\Z(0)$: this proves (3) and the fact that 
each $H_i$ is connected by Hopf's theorem. Moreover from the inclusion 
$L\subset K\subseteq H$, we
get a fibration  $$S^1\hookrightarrow H/L\to H/K.$$
Since $H/L$ is a sphere (see e.g. [AA]), we get that $H/K$ must be diffeomorphic
to 
some complex projective
space and this proves (4).\par 
We now want to prove that $H_1\cap H_2=K$. Indeed, 
we first note that $H_1,H_2,K$ share a common 
maximal torus $T$ of $G$ and $Z(H_i)\subseteq Z(K)\subset T$; 
now it is easy to see that 
$H_1\cap H_2$ is the centralizer in $G$ of the torus 
$Z(H_1)\cdot Z(H_2)\subset T$, 
hence it is 
connected with Lie algebra $\h_1\cap \h_2$. If $K\not= H_1\cap H_2$, 
we may find a non 
zero vector $v\in\m\cap \h_1\cap\h_2$; for such $v$, 
we consider  the function 
$\varphi(t)=||\hat v||^2_{\gamma_t}$ for $t\in[0,\d]$. Since 
$\varphi(0)=\varphi(\d)=0$, by 
Rolle theorem we get that $\varphi'(c)=0$ for some 
$c\in(0,\d)$; but, by 
$$\varphi'(t) = -a(t) \B([Z, Jv], v) = \pm ||\hat Z||_{\gamma_t}  \B([Z, Jv], 
v)$$
$\varphi'(t)$ 
cannot vanish in the open interval $]0,\d[$ and we are done.\par
We now come to our last claim about geodesic symmetries. It is known (see 
[AA]) that at any singular point $\gamma_{k\d}$ we can find an element $\sigma$
belonging to the stabilizer $G_{\gamma_{k\d}}$ such that $\sigma$ reverses the 
normal geodesic $\gamma$; it is clear that $\sigma$ normalizes $L$. In our case
we have shown that each singular orbit is a complex submanifold; hence the slice 
representation of, say, $H_1$ must preserve some complex structure on the normal
space.
It then follows that, if $\nu_1$ denotes the slice representation, then 
$\nu_1(H_1)$ is one of the following groups: $SU(m),U(m),Sp(m),T^1\cdot Sp(m)$.
Now,
in each of these cases, the geodesic symmetry lies in the connected component
$N_H^o(L)$ and 
therefore $\sigma\in N^o_G(L)=K$.\qed\enddemo
\bigskip
\remark {Remark 4.2} The fact that the singular orbits are 
complex submanifolds relies on the
projectability of the complex structure $J$. Indeed the cohomogeneity 
one action of $G=SO(n+1)$ 
on the complex quadric $Q_n$ admits a totally real singular orbit. 
Note that a singular orbit could
be neither complex nor totally real; such an example is given by a 
cohomogeneity one action of the 
group $G=U(n+1)$ on $Q_{2n}$ (see [CN]).\par
It is also known (see e.g. [HS]) that, given a 
cohomogeneity one K\"ahler $G$-manifold, we can always blow up the singular
orbits 
and reduce ourselves to the case when both singular orbits have complex
codimension one.
However, as far as we know, there is no control on the differential-geometric
properties 
of the metric when we blow up. At this regard, we refer the reader also to the 
paper by Koiso-Sakane ([KS]).\endremark
\bigskip
\bigskip
We want now to determine
a fundamental relation between the invariants $\Z$ and $a$. We begin
with a technical lemma.\par 
\proclaim{Lemma 4.3} Let $\gamma$ be parametrized
so that (2.1) holds and let $a$ be defined as in (3.13).
Then (up to a reparametrization $t\to -t$)
 for any $t\in A$, 
$$\frac{1}{a(t)} \hat Z_{\gamma_t} = J\gamma'_t\tag 4.4$$
\endproclaim
\demo{Proof} First of all, observe that both 
$\frac{1}{a(t)} \hat Z_{\gamma_t}$
and $J\gamma'_t$ are ad($\l$)-invariant vectors
in $T_{\gamma_t} (G\cdot \gamma_t)$. By means of the identification
(3.2), they have to correspond to two ad($\l$)-invariant vectors
in $\n$ and hence, since the action of $G$ is ordinary, they both 
belong to $\R Z$. Furthermore, $|a(t)| = ||\hat Z||_{\gamma_t}$
and hence   both vectors are of unit length. So
we may assume that (4.4) holds
 for all $t\in
]-\d, 0[$. On the other hand, for $t\in ]0, \d[$, we may consider
the geodesic symmetry $\sigma\in H_1$, which lies in $K$ by Prop. 4.1. So, 
$$ J\gamma'_t = - \sigma_* (J \gamma'(-t)) =
 \frac{1}{||\hat Z||_{\gamma(-t)}}
\sigma_* (\hat Z_{\gamma(-t)}) = 
\frac{1}{||\hat Z||_{\gamma_t}}\widehat{Ad(\sigma) Z}_{\gamma_t} =
\frac{1}{||\hat Z||_{\gamma_t}}\hat Z_{\gamma_t}$$
where we have used the fact that $Z \in \z(\k)$. 
This proves (4.4) in the interval $]0, \d[$. By reflecting
the geodesic at all singular points, the same argument proves
(4.4) at all $t\in A$.\qed
\enddemo
\bigskip
We continue with an important proposition which will bring to a complete
description of the invariant $\Z$.\par
\proclaim{Proposition 4.4} Let $\gamma$ parametrized so that
 $\gamma_0$ and $\gamma_\d$ are the intersections
with the two singular orbits. Then
$\Z(t)$ is even with respect to the reflections
of the parameter $t\to -t$ and $t\to 2\d - t$. In particular, $\Z(t)$ and
$a(t)$ are both periodic of period $T = 2\d$.
\endproclaim
\demo{Proof} 
To see that 
$\Z$ is invariant with respect to the reflection $t\to -t$
(or $t\to 2\d - t$), consider the geodesic symmetry $\sigma$  
at the singular point $\gamma_0$ (or $\gamma_\d$)
such that 
$\gamma(-t) = \sigma\cdot \gamma_t$.
We know that 
$\sigma\in K$ by Prop. 4.1. Moreover,  using the definition of
$\Z$ we obtain that
$$\Z(-t) = Ad(\sigma)\Z(t) = \Z(t)$$
because $\Z(t) \in \z(\l) + \a$. This proves the claim.\qed
\enddemo
We conclude by summing up the results of Propositions 4.1, 4.2 and 4.3 
in order to  get the
following straightforward
characterization of the image $\Z(\R)$. For the statement,
we need to  recall a few basic fact on the flag 
manifold $G/K$ (see e.g. [Be]). 
An element $X\in \z(\k)$ is called {\it regular\/}
if $C_G(X) = K$. The set of regular elements is open and 
dense in $\z(\k)$ and each connected component is called
{\it $T$-chamber\/} of the flag manifold. Let $\t$ be a maximal torus
for $\k$ (and hence also for $\g$) and $R$ (resp. $R_K$) be the root system
of the pair $(\g^\C, \t^\C)$ (resp. of $(\k^\C, \t^\C)$). Then
$X\in \z(\k)$ is regular if and only if 
$\alpha(X) \neq 0$ for all $\alpha \in R\setminus R_K$. Finally,
recall that
for any $G$-invariant complex structure $J$ on $G/K$ there exists
an ordering of the roots so that, for any $\alpha \in R\setminus
R_K$,  $J E_\alpha = +i E_\alpha$
if and only if $\alpha\in R^+$. The $T$-chamber of the elements
$\{ iX\in \z(\k)\ : \alpha(X) >0\}$ is uniquely associated
to $J$ and it is called the {\it positive $T$-chamber
corresponding to $J$\/}.\par
\bigskip
\proclaim{Corollary 4.5} Let $G/K$
be the flag manifold which is the image by the moment
map of the regular orbit $G/L$ of $M$.
Then the trace $\ell_Z = \Z(\R)$ is a segment parallel to $Z$, which 
verifies the following properties:
\roster
\item"i)" $\ell_Z$ is contained in the 
closure of the positive $T$-chamber corresponding to $J$;
\item"ii)" the centralizers $H_1$ and 
$H_2$ of the endpoints $Z_1$ and $Z_2$ of 
$\ell_Z$ 
are such that $H_1\cap H_2 = K$ and $H_i/K\simeq \C P^{m_i}$
for some $m_i$ (in case $Z_i$ is inner to the $T$-chamber,
$H_i/K$ is point and we consider $m_i = 0$);
\item"iii)" there exist two invariant complex structures $J_i$, $i = 1,2$
on the flag manifolds $G/H_i$ such that the projections
$$\pi_1\: (G/K,J_K)\to(G/H_i,J_i)$$
 are holomorphic; 
\item"iv)" the length of $\ell_Z$ is equal to $\int_0^\d||\hat Z||_{\gamma_t}
dt$.
\endroster
\endproclaim
\demo{Proof} The only claim that it is not an immediate
corollary of the previous results is the fact that $\ell_Z$ lies 
in the positive $T$-chamber of $J$; but this is a consequence
of the fact that
the symmetric two form on $\m$
$$g_t(X,Y) \= \B([\Z(t), J_K X], Y)$$
is positive definite  at all $t\in \R$.\qed
\enddemo
\bigskip
 \bigskip
\subhead 5. Abstract models: proof of the 
Theorem 1.2
\endsubhead
\bigskip
In the following,  $G$ is a compact semisimple Lie group,
$G/K$ is a flag manifold and $L$ is a closed
subgroup of codimension one in $K$.
 As usual we will consider the 
decomposition 
$$\g = \l + \n = \l + \R Z + \m$$
with the  meaning of the symbols as in the previous sections.
Assume also that the  pair $(G,L)$ is ordinary. We may now
give the definition of abstract model:\par
\bigskip
\definition{Definition 5.1} Let $J_K$ be a complex
structure on $G/K$ and let $\ell_Z$ an
oriented segment in $\z(\k)$
which is parallel to $Z$. We say that it is an {\it admissible
segment for $(G,L, J_K)$\/} if it verifies i), ii) and iii) of Corollary 
4.5.\par
Denote by $Z_i$, $i=1,2$,  the endpoints of $\ell_Z$:
we call {\it degree\/} of $Z_i$ the integer $\deg(Z_i) = m_i+1$, where 
$m_i$ is the complex dimension of  
$H_i/K \simeq \C P^{m_i}$ (note: in case $Z_i$ is inner
to the $T$-chamber, we set $\deg(Z_i) = 1$).\par
If $C = ||Z_2 - Z_1||$, w.r.t. $\B$, 
a smooth function $f\: \R \to [0, C]$ is said {\it admissible 
parametrization
of $\ell_Z$\/}
if:  
\roster
\item"a)" $f(0) = 0$, $f(\d) = C$
and it is monotone increasing in the interval
$]0, \d[$, for some $\d$; 
\item"b)" it is invariant by the
symmetries $t\to -t$ and $t \to 2\d - t$; 
\item"c)" $f''(0) = 1 = - f''(\d)$.
\endroster
A pentuple of the form $\Cal K = (G,L, J_K, \ell_Z, f)$
where $\ell_Z$ is an admissible segment and $f$ is an admissible
parametrization of $\ell_Z$ is called
{\it abstract model  for a compact K.O.P. manifold\/} (=
K\"ahler manifold with Ordinary and Projectable
 cohomogeneity one action).\par
We say that two abstract models $ (G, L, J_K, \ell_Z, f)$
and $ (G', L', J_{K'}, \ell_{Z'}, f')$ are {\it equivalent\/}
if and only if: i) $G= G'$ and $(G/K, J_K)$ and $(G'/K', J_{K'})$ are 
biholomorphic
via an automorphism  $\varphi$ of $G$; ii) $\varphi(L) = L'$; iii) 
$\varphi_*(Z) = Z'$ and $\varphi_*(\ell_Z) = \ell_{Z'}$ with 
the same orientation, or $\varphi_*(Z) = - Z'$ and $\varphi_*(\ell_Z)$
coincides with $\ell_{Z'}$ but with 
opposite orientation; iv) $f = f'$.
\enddefinition
\bigskip
The proof of Theorem 1.2 coincides 
with the  proof of following Theorem 5.2, since the second claim
of Theorem 1.2 is a direct consequence of the definitions.\par
\medskip
\proclaim{Theorem 5.2} For any abstract model 
$\Cal K = (G, L, J_K, \ell_Z,f)$, there exists
a compact K.O.P.  $G$-manifold $(M,g,J)$ such that:
\roster
\item"a)" the invariant associated to the complex structure is 
${\Cal J}(t) = (Z, J_t, Z)$,
 where 
$J_t$ coincide for all $t\in A$ with the unique  
complex structure which projects onto the 
complex structure $J_K$ of $G/K$;
\item"b)" the invariant associated to the K\"ahler form is $\Z(t)
= Z_1 - f(t)Z$, where $Z_1$ is the first endpoint of $\ell_Z$;
\item"c)" the invariant 
ociated to the metric tensor is the function $a(t)
= f'(t)$
\endroster
We call such a $(M,g,J)$ a {\it realization\/} of the abstract model $\Cal 
K$.\par
Any compact K.O.P. manifold $(M,g,J)$ is
$G$-equivalent to a realization of an abstract model $\Cal K$
and this corresponding
abstract model  is the same (up to equivalence of 
abstract models)  for any other
manifold in the same $G$-equivalence class.
\endproclaim
\demo {Proof} We only prove the
existence of a realization for any abstract model,
since the proof of 
uniqueness up $G$-equivalence is straightforward.
As usual, we will denote by $H_i = C_G(Z_i)$ the centralizers of the 
endpoints of $\ell_Z$.\par
By definitions, the triple ($H_1,L,H_2$) defines a 
$G$-manifold which is of cohomogeneity one 
w.r.t. $G$, with regular orbits equivalent
to $G/L$ and with two singular orbits $G/H_1$ and
$G/H_2$ respectively (see [AA] for a detailed exposition). 
More precisely, each subgroup $H_i$, $i=1,2$, has an
orthogonal representation $\rho_i:H_i\to O(V_i)$ in some vector 
space $V_i$ such that $\rho_i(H_i)$
acts transitively on the unit sphere $S(V_i)$ and $S(V_i)\cong H_i/L$; 
the manifold $M$ is 
then obtained by glueing together two disk bundles over the two 
singular orbits $G/H_i$. Since
the two singular orbits are flag manifolds, hence simply connected, 
by Seifert-Van Kampen
Theorem, the resulting manifold $M$ will be simply connected.\par 
We begin with a Lemma which will be
useful later on. \medskip \proclaim {Lemma 5.3} For $i=1,2$, 
we have that $\rho_i(H_i)$ is
$SU(k_i),U(k_i)$ or $T^1\cdot Sp(k_i)$ for some $k_i\geq 1$.\endproclaim 
\demo{Proof} The proof is obtained by inspection of Borel's 
list of compact Lie groups $H$ acting 
transitively on some sphere with isotropy $L$, when we consider 
only ordinary pairs ($H,L$).\qed
\enddemo
\bigskip
We further recall that each singular orbit $G/H_i$, $i=1,2$ 
has a $G$-invariant tubular neighborhood 
which is $G$-diffeomorphic to the total space of the vector bundle 
$G\times_{H_i}V_i$ over
$G/H_i$.\par
We may identify the regular part $M_{\text{reg}}$ with $G/L\times ]0,\d[$ 
and we shall consider the 
curve $\gamma:]0,\d[\to M_{\text{reg}}$ given by $\gamma(s)=(eL,s)$.
 Given the splitting 
$$\g = \l +\Bbb R Z + \m,$$
we can identify the tangent space $T_pM_{\text{reg}}$ with 
$\Bbb R{{\partial}\over{\partial
s}}+\Bbb R Z + \m$; moreover we choose a basis $\{e_1,\dots,e_m\}$ of 
$\m$ so that 
$\{{{\partial}\over{\partial s}},\hat Z,\hat e_j\}_{\gamma(s)}$ is a 
frame along $\gamma$. We shall 
denote by $ds$ and $\eta$ the 1-forms which are dual to 
${{\partial}\over{\partial s}},\hat Z$
respectively. \par
We may now define a $G$-invariant metric $g$ on $M_{\text{reg}}$ by 
$$g = ds^2 + f'(s)^2\eta^2 + g|_{\hat\m\times\hat\m},\tag 5.1$$
where, for $X,Y\in\m$,
$$g|_{\hat\m\times\hat\m}(\hat X,\hat Y) = \B([\Z(s),J_KX],Y),\tag 5.2$$
where $\Z(s) = Z_1-f(s)Z$.
We now prove that $g$ extends smoothly to a $G$-invariant 
Riemannian metric on the whole $M$. In
order to do that, we restrict $g$ to $M_1=G\times_{H_1}V_1$.\par
First of all, we note that $\m$ admits a further splitting as 
$\m=\m_1+\tilde\m$, where 
$\h_1=\k+\m_1$; the summand $\tilde\m$ defines a smooth distribution 
$\Cal H$ on the whole $M_1$, so
that (5.2) and property (i) in Cor. 4.5 allow to extend smoothly 
$g|_{\Cal H\times\Cal H}$.\par
We are therefore left with proving that the restriction of $g$ on the slice 
$\bV$
in $M_1$ defined as
$$\bV = \{ [(e; v)]\ v\in V_1\ \}\subset G\times_{H_1}V_1, $$
extends smoothly on the whole $\bV$.\par
We denote by $g_o$ the standard Euclidean metric on $V_1$, 
which is invariant by the linear
group $\rho_1(H_1)$ and we express it as
$$g_o = ds^2 + s^2\eta^2 + g_o|_{\hat\m_1\times\hat\m_1}.$$
We observe that both $g|_{\hat\m_1\times\hat\m_1}$ and 
$g_o|_{\hat\m_1\times\hat\m_1}$ are 
Ad($K$)-invariant; according to Lemma 5.3, the Ad($K$)-module $\m_1$ is 
either irreducible or 
it splits as the sum of two irreducible, 
inequivalent submodules $\m_1^{(0)},\m_1^{(1)}$. Therefore
we my write that
$$g|_{\hat\m_1\times\hat\m_1} = \lambda^2(s) 
g_o|_{\hat\m_1\times\hat\m_1}\ ,$$
if $\m_1$ is Ad($K$)-irreducible, or otherwise
$$g|_{\hat\m_1\times\hat\m_1} = 
\lambda^2(s) g_o|_{\widehat{\m_1^{(0)}}\times\widehat{\m_1^{(0)}}}+
\mu^2(s) g_o|_{\widehat{\m_1^{(1)}}\times\widehat{\m_1^{(1)}}}$$
for suitable positive functions $\lambda,\mu$. 
From (5.2) and the properties of $\Z=Z_1-fZ$, it is clear
that  $\lambda(s)$ and $\mu(s)$ admit smooth and even 
$C^\infty$-extensions over $\Bbb R$.\par
Using the results in [Ve], we have that $g|_{T\bV}$ extends 
smoothly as a Riemannian metric on the
whole $T\bV$ if and only if the functions 
${{f'(s)}\over{s}},\lambda(s),\mu(s)$ extend to
smooth, even and positive functions over $\Bbb R$ with 
$\lim_{s\to 0}{{f'(s)}\over s} = 1$ and these conditions are satisfied 
in view of the previous remarks and our
hypothesis on $f(s)$.\par
In the same way we prove that the metric 
$g$ extends smoothly on the whole $M_2=G\times_{H_2}V_2$.\par
We now want to define an almost complex structure $\bJ$ 
on $M$. We first define it on the regular
part by putting, along the curve $\gamma$,
$$\bJ({{\partial}\over{\partial s}}) = {1\over{f'(s)}}\hat Z,\quad 
\bJ(\hat Z) = -f'(s){{\partial}\over{\partial s}},\quad
\bJ|_{\hat m} = J_K|_{\hat m}\ .$$
It is clear that the metric $g$ is $\bJ$-Hermitian. 
We want to prove that $\bJ$ extends to a smooth 
almost complex structure on the whole $M$. 
Using property (iii) in Cor. 4.5 and the same arguments
as in the proof of extendability of the metric, 
we are left with proving that $\bJ$ extends smoothly
on the whole slice $\bV$ in $M_1$. The restriction of 
$\bJ$ on $T\bV$ takes the form
$\bJ({{\partial}\over{\partial s}}) = {1\over{f'(s)}}\hat Z$, 
$\bJ(\hat Z) = -f'(s){{\partial}\over{\partial s}}$,
$\bJ|_{\hat m_{1}} = J_K|_{\hat m_{1}}$.
Following Sakane ([Sa]), we change the parameter 
$s$ by considering a new function 
$$r(s) = \exp(\int_{s_o}^s{1\over{f'(u)}}du), \quad s\in ]0,\d[$$
where $s_o\in ]0,\d[$ is a fixed value of the parameter.\par
We have 
\medskip
\proclaim {Lemma 5.4} The function $r$ admits a 
smooth and odd extension to $]-\d,\d[$. 
\endproclaim
\demo{Proof} Indeed, by assumption, we know that the function 
$f'(u)$ admits an odd extension on
the whole $\Bbb R$, so that we may write $f'(u) = u(1+u^2b(u))$, 
for some $C^\infty$, even function
$b$. Now, for $s\in ]0,\d[$ we have 
$$r(s) = {s\over{s_o}}\exp(-\int_{s_o}^s {{ub(u)}\over{1+u^2b(u)}}du)\tag 5.3$$
Since $f'$ does not vanish on $]-\d,0[$, we have that the integral 
$\int_{s_o}^s
{{ub(u)}\over{1+u^2b(u)}}du$ has a natural $C^\infty$ extension 
on the whole $]-\d,\d[$ as an even
function. Then (5.3) shows that $r$ extends as claimed.\qed\enddemo
\bigskip
If we now define a map $\phi$ on 
$U^* = \{v\in \bV\setminus \{[e,0]\};\ ||v||<\d\}$, where 
the norm $||\cdot||$ refers to the Euclidean metric $g_o$, by
$$\phi(v) = {{r(||v||)}\over{||v||}}v,$$
then the property of the function $r$ established in 
Lemma 5.4, shows that $\phi$ extends to a
diffeomorphism of $U=U^*\cup \{[e,0]\}$ onto some 
symmetric neighborhood of $0$ in $\bV$.
It is now
easy to verify that  
$$\eqalign{(\phi_*\bJ)({{\partial}\over{\partial s}}) &= 
{1\over{s}}\hat Z\cr
(\phi_*\bJ)|_{\hat m_1} &= J_K|_{\hat\m_1}.\cr}$$
It then follows that the almost complex structure 
$\phi_*\bJ$ can be extended to the whole 
$T\bV$, since it coincides with the unique 
(up to sign) complex structure on $V_1$ which is 
invariant by the linear group $\rho_1(H_1)$ (see Lemma 5.3).\par
Again, the same arguments can be applied to 
the submanifold $M_2=G\times_{H_2}V_2$, proving that 
$\bJ$ extends smoothly to an almost complex structure on $M$. 
We now prove that $\bJ$ is integrable.
\bigskip
\proclaim {Lemma 5.5} The almost complex structure $\bJ$ 
is integrable\endproclaim
\demo{Proof} We will check this, 
proving that the Nijenhuis tensor $N$ of $\bJ$ vanishes 
on the regular part, which is
dense in $M$. \par
The regular part $M_{\text{reg}}$ fibers as
$$\pi:M_{\text{reg}}\cong G/L\times ]0,\d[\to G/K,$$
with typical fibre $F\cong S^1\times ]0,\d[$; since the vertical 
space $\Cal V$ of $\pi$ is $\bJ_1$-invariant and the fibres are 2-dimensional, 
we have that $N(\Cal V,\Cal V)=0$. \par
Moreover, the equation $N(\hat\m,\hat\m)=0$ is automatic from the 
fact that $\bJ|_{\hat\m}$ projects down to the integrable complex 
structure $J_K$ on the flag manifold $G/K$. 
We are therefore left with proving that $N(\Cal V,\hat\m)=0$.\par
If we denote by $\xi$ the vector field ${{\partial}\over{\partial 
s}}$, where $s\in ]0,\d[$, then it is enough to prove that 
$$N(\xi,\hat\m)=0 \tag 5.4$$
We choose a basis $\{e_1,\dots,e_m\}$ of $\m$ and compute 
$$\eqalign {N(\xi,\hat e_{i}) &= [\bJ\xi,\bJ\hat e_{i}] - 
[\xi,\hat e_{i}] - 
\bJ[\xi,\bJ\hat e_{i}] - \bJ[\bJ\xi,\hat e_{i}]\cr 
{}& = [\bJ\xi,\bJ\hat e_{i}] - \bJ[\xi,\bJ\hat e_{i}]\ .
\cr}$$
Now, we fix a regular point $p\in M_1{}_{\text{reg}}$ with $G_p=L$ 
and we write locally around $p$ 
$$\bJ\hat e_{i} = a\xi + b\hat Z + \sum_{k=1}^m 
c_{k}\hat e_{k},$$
for some $C^\infty$ functions $a,b,c_{k}$; from the construction of 
$\bJ$ it follows that, if $\gamma$ 
denotes the integral curve of $\xi$ through $p$, then 
$a\circ\gamma = b\circ \gamma =0 $ and the $c_{k}\circ\gamma$ are 
constant functions for $k=1,\dots,m$.\par
It is then straightforward to see that $[\xi,\bJ\hat e_{i}]_{p}=0$. 
Therefore, 
$$N(\xi,\hat e_{i})_{p} = [\bJ\xi,\bJ\hat e_{i}]_{p}.$$
If we now write locally 
$$\bJ\xi = \beta\hat Z + \sum_{k=1}^m \mu_{k}\hat e_{k}$$
for some $C^\infty$ functions $\beta,\mu_{k}$ with $\mu_{k}(p)=0$, we 
can compute
$$\eqalign{[\beta\hat Z + \sum_{k=1}^m \mu_{k}\hat e_{k},\bJ\hat 
e_{i}]_{p}&= (-(\bJ\hat e_{i})(\beta)\hat Z + \beta [\hat Z,\bJ\hat e_{i}] -
(\bJ\hat e_{i})(\mu_{k})\hat e_{k})_{p} = \cr
{}&= (-(\widehat{J_K e_{i}})(\beta)\hat Z + \beta[\hat Z,\bJ\hat e_{i}] -
(\widehat{J_K e_{i}})(\mu_{k})\hat e_{k})_{p}\cr}\tag 5.5$$
We now prove the following 
\proclaim {Sublemma 5.6} We have 
$$[\hat Z,\bJ\hat e_{i}]_{p}= [\hat Z,\widehat{J_Ke_{i}}]_{p}.$$ 
\endproclaim
\demo{Proof} Indeed
$$\eqalign{ [\hat Z,\bJ\hat e_{i}]_{p}&= \bJ[\hat Z,\hat 
e_{i}]_{p} = -\bJ\widehat{[Z,e_{i}]}_{p} = \cr
{}& = -\widehat{J_K[Z,e_{i}]}_{p} = -\widehat{[Z,J_Ke_{i}]}_{p} = [\hat 
Z,\widehat{J_Ke_{i}}]_{p},\cr}$$
where we have used the fact that $J_K$ is ad($Z$)-invariant.\qed\enddemo
\enddemo
From Sublemma 5.6 and (5.5), we conclude that
$$N(\xi,\hat e_{i})_p = [\bJ\xi, \bJ\hat e_i]_p = [\bJ\xi, \widehat{J_K e_i}]_p
= 0\ .\qquad\qquad \qed$$
\medskip
We are left with proving that the metric $g$ is K\"ahler, that is the associated
fundamental form $\omega$ is closed. Again, we will check 
this on the regular part. 
The condition $d\omega=0$ on 
$M_{\text{reg}}$ is equivalent to the following four equations along 
the normal geodesic
$\gamma(s)$:\par \roster
\item"i)"\ $d\omega(\hat v_1,\hat v_2,\hat v_3)=0$;
\item"ii)"\ $d\omega(\xi,\hat v_1,\hat v_2)=0$; 
\item"iii)"\ $d\omega(\hat Z,\hat v_1,\hat v_2)=0$; 
\item"iv)"\ $d\omega(\xi,\hat Z,\hat v_1)=0$,
\endroster
where $v_1,v_2,v_3$ belong to $\m$ and $\xi = {{\partial}\over{\partial 
s}}$.\par
We note that equations (i) and (iii) are equivalent to condition (3.11).
Equation (ii) is equivalent to 
the equation (4.1) satisfied by the invariant $\Z(s)$.\par
As for the last equation (iv), using 
again the fact that ${\Cal L}_{\hat v_1}\omega = {\Cal L}_{\hat Z}\omega = 0$, 
we get
$$\eqalign { d\omega(\xi,\hat Z,\hat v_1) &= \xi\omega(\hat Z,\hat v_1) - 
\hat Z\omega
(\xi,\hat v_1) + \hat v_1\omega(\xi,\hat Z)\cr
{}&-\omega([\xi,\hat Z],\hat v_1) + \omega([\xi,\hat v_1],\hat Z) - 
\omega([\hat Z,\hat
v_1],\xi)=\cr
{}& = \xi\omega(\hat Z,\hat v_1) + \omega(\xi,[\hat v_1,\hat Z]) =\cr
{}& = {d\over{ds}}\B([\Z(s),Z],v_1) + g_s(Z,[v_1,Z]) = 0,\cr}$$
since $\Z(s)$ and $Z$ commute and $Z$ is $g_s$-orthogonal to $\m$. \par
This concludes the proof of the theorem.\qed\enddemo
\bigskip
\bigskip
\subhead 6. The Ricci tensor of a realization of $\Cal K = (G, L, 
J_{K},\ell_{Z},f)$
\endsubhead
\bigskip
In this section, we will compute the Ricci tensor of the Riemannian manifold
($M,g$) corresponding to an abstract model $\Cal K = (G, L,
J_{K},\ell_{Z},f)$.\par
We shall keep the same notations as in the previous section. Moreover, we 
denote by $Z_{1},Z_{2}$ the endpoints of $\ell_{Z}$, so that $\Z(t) = 
Z_{1}-f(t) Z$. 
Also, we fix a basis $\{ v_1, \dots, v_p\}$ for $\m$ which is orthonormal 
with respect to $\B$. If $g$ is the Riemannian metric
determined by $\Z$, the norms $||\hat v_i||$ can be determined
by the formula
$$||\hat v_i||^2 = 
\B([Z_1, J v_i], v_i) - f(t) \B([Z, J v_i], v_i) \=
h_i - f(t) k_i \ .$$
Along the geodesic $\gamma(t)$, at all the regular points, 
we will denote by $e_i$
the normalized vectors parallel to the $\hat v_i|_t$, that is
$$e_i = \frac{1}{||\hat v_i||} \hat v_i =
\frac{1}{\sqrt{h_i - f k_i}} \hat v_i\ .$$
We recall that $||\hat Z||^2 = a(t)^2 = f'(t)^2$. 
Finally, recall also that the Riemannian metric $g$ verifies 
$g(\hat Z, \hat v_i) \equiv 0$
for all $v_i$. We will also denote by $g_t$ the Riemannian metric induced
on the regular orbit $G\cdot \gamma_t$.\par
\bigskip
Let us now start the computation of the Ricci curvature at the regular points.
First of all we need the following Lemma (see also [Sa]):\par
\bigskip
\proclaim {Lemma 6.1} Let $r$ be the Ricci tensor of the 
$G$-invariant metric
$g$ on a realization $M$ of $\Cal K$. Then $r(\xi,\hat\m)|_{\gamma_{t}}=0$.
\endproclaim
\demo{Proof} We recall that $[\l,\m]=\m$. 
If $v\in\m$, then $v=\sum_i[X_i,Y_i]$
for some $X_i\in\l$ and $Y_i\in\m$. Then
$$r(\xi,[\hat X_i,\hat Y_i])|_{\gamma_{t}} = 
\hat X_i r(\xi,\hat Y_i)|_{\gamma_{t}} -
 r([\hat X_i,\xi],\hat Y_i)|_{\gamma_{t}} = 0$$
since $X_i\in\l$. \qed \enddemo
\bigskip
Lemma 6.1 is not the only vanishing 
property of $r$. 
Indeed, note that the bilinear forms induced
by  $g_t$ and  $r$ on $\m$,
 via the identification
map $\phi_t$ in (3.2), are both $ad(\k)$-invariant.
Now, if $\a_t$ is an $ad(\k)$-invariant, symmetric, bilinear form
on $\m$ and if $\m=\m_1+\dots+\m_p$ is the 
ad($\k$)-invariant decomposition into
irreducible submodules, 
we have that $\a_t(\m_i,\m_j) \equiv 0$ for $i\not=j$ and that
$\a_t|_{\m_i}$ is a
 multiple of $\B|_{\m_i}$ for all $i=1,\dots,p$.\par
From such observations, it follows immediately that
$$g_t(\hat v_i, \hat v_j) = 0\quad,\quad
 r(\hat v_i, \hat v_j)|_{\gamma_{t}} = 0 \qquad \forall i \neq j.$$
In conclusion, $r|_{\gamma_t}$ is uniquely
determined by the values of  $r(e_i, e_i)$ and $r(\xi,\xi)$.\par 
We may now compute the Ricci tensor by considering the
Riemannian submersions 
$\pi\: (M_{\text{reg}},g) \to (\R, dt)$ and 
$\kappa\: (G\cdot\gamma(t), g_t) \to (G/K, g^{\kappa}_t)$, where 
$g^{\kappa}_{t}$ is the $G$-invariant metric on $G/K$ induced by 
$g_{t}$. The reader 
interested in detailed computations is referred to [PS].\par
We will also denote by $A^\pi$, $A^\kappa$, $T^\pi$, $T^\kappa$
\dots the usual
O'Neill tensors (see e.g. [Be]),
which are related to the submersions $\pi$ and
$\kappa$, respectively. \par
By O'Neill's formulae for the submersion $\kappa$, we 
compute $r_t(e_i, e_i)$:
$$r_t(e_i, e_i) = r^\kappa_t(e_i, e_i) - 2 g_t(A^\kappa_{e_i} J\xi,
A^\kappa_{e_i} J\xi) =$$
$$ = r^\kappa_t(e_i, e_i) - 
\frac{(f')^2}{2} \frac{k_i^2}{(h_i - f k_i)^2}$$
Now, in order to compute $r(e_i, e_i)$ we  use  O'Neill's formulae
for the submersion $\pi$:\par
$$r(e_i, e_i) = 
r_t(e_i, e_i) -
g(N^\pi, T^\pi_{e_i} e_i) + (\tilde \delta T^\pi)(e_i, e_i).$$
A lenghty but straightforward computation shows that 
$$r(e_i, e_i)  = r^\kappa_t(e_i, e_i) + 
\left[
f''  -
\frac{(f')^2}{4} 
\sum_m\frac{k_m}{h_m - f k_m}\right] 
\frac{\B([Z, J v_i], v_i)}{||\hat v_i||^2},$$
so that for all $X,Y\in \m$
$$\boxed{
r(\hat X, \hat Y) =
r^\kappa_t(\hat X, \hat Y) + \left[
f''  -
\frac{(f')^2}{4} 
\sum_i\frac{k_i}{h_i - f k_i}\right] 
\B([Z, J X], Y)}\tag6.1
$$
\bigskip
The evaluation of $r(\xi, \xi)$ can also be computed by using O'Neill's
formulas for the submersion $\pi$. So we get that 
$$ \boxed{r(\xi, \xi) =
- \frac{f'''}{f'}
+ \frac{1}{2} f'' \sum_i \frac{k_i}{h_i - f k_i}
+ \frac{1}{4} (f')^2 \sum_i \frac{k_i^2}{(h_i - f k_i)^2}}\tag6.2 $$
\medskip
$$ \left(= - \frac{1}{f'}
\frac{d}{dt}\left[f'' - 
\frac{(f')^2}{4} \sum_i \frac{k_i}{h_i - f k_i}\right]
\right)\tag6.$2_1$ $$
We remark here that, using the fact that the 2-form $\rho$ is closed and Lemma
6.1, if $X$ and $Y$
belong to $\m$, then
$$\xi\rho(\hat X,\hat Y) = -g([\hat X,\hat Y],J\xi)r(\xi,\xi).\tag 6.3$$
It then follows that, if the Einstein equation $r(\hat X,\hat Y) = cg(\hat X,\hat
Y)$ is 
fulfilled for some $c\in \Bbb R^+$ and for all $X,Y\in \m$, then also 
the Einstein equation 
$r(\xi,\xi) = c$ is satisfied, by (6.3) and (4.1).
\bigskip\bigskip

\subhead 7. Einstein-K\"ahler manifolds of cohomogeneity one
\endsubhead
\bigskip
Recall that if $(G/K, J_K)$ is a flag manifold with invariant complex
structure $J_K$, then for any invariant K\"ahler matric $g^\kappa$ compatible 
with
$J_K$, the Ricci form $\rho^\kappa$ can be written as 
$$\rho^\kappa(\hat X, \hat Y)_{eK} =
\B([Z^\kappa, J_K X], Y),\quad X,Y\in \m\tag 7.1$$
where $Z^\kappa \in \z(\k)$ {\it and it does not 
depend on the metric $g$ but only on the complex structure
$J_K$\/} (see e.g. [BFR]). We will call it the {\it Ricci 
invariant of $(G/K, J_K)$\/} (the explicit expression of 
$Z^{{\kappa}}$ is given in the Introduction).\par
Using this fact and 
 (6.1) and (6.3) we have that 
a realization of $\Cal K$ is a K\"ahler-Einstein manifold with Einstein
constant c (i.e. $r = c g$) if and only if the invariant
$\Z(t) = Z_1 - f(t) Z$
verifies the following ordinary differential equation
$$
\left[
f''  -
\frac{(f')^2}{4} 
\sum_i\frac{k_i}{h_i - f k_i} + cf\right] Z =   cZ_1 -
Z^\kappa \tag7.2$$
\medskip
We will now look for admissible segments and parametrizations which
solve (7.2) for $c=1$.\par
First of all, there is a constant $D$
such that
$$  Z_1 - Z^\kappa = - D Z \tag7.3$$
$$f''  -
\frac{(f')^2}{4} 
\sum_i\frac{k_i}{h_i - f k_i} + f + D= 0 \ .\tag7.4$$
Let us call $S_1$ and $S_2$ the two singular orbits of a
realization of $\Cal K$ and recall that $S_i = G/H_i$,
where $H_i = C_G(Z_i)$. Note also that the complex codimension
of each $S_i$ coincide with $m_i=\deg(Z_{i})$ (see Def. 5.1 and 
Cor. 4.5). Furthermore,
this implies that $2(m_1-1)$ equals the number 
of constants $h_i$ which vanish and that $2(m_2-1)$ is the number of 
vanishing values $h_i - f(\d) k_i$, where $\d$ is the distance
between the two singular orbits.\par
  Now if $f$ is a solution
of (7.4) for some open interval $]0, t_o[$, then it represents an 
admissible parametrization only if 
$$\lim_{t\to 0^+}f(t) = 0 =  \lim_{t\to 0^+}f'(t)\quad
\lim_{t\to 0^+}f''(t) = 1$$
and this implies that 
$$D = -(\lim_{t\to 0^+}f''(t) - 
\lim_{t\to 0^+}\frac{(f')^2}{4} 
\sum_{i=1}\frac{k_i}{h_i - f k_i}) = - m_1\tag7.5$$
The analogous conditions for a solution defined in an open interval
$]t_o, \d[$ are
$$\lim_{t\to \d^-}f(t) = f_\d>0\quad
\lim_{t\to \d^-}f'(t) = 0\quad
\lim_{t\to \d^-}f''(t) = -1$$
and from this we get that
$$f_\d = -(\lim_{t\to \d^-}f''(t) - 
\lim_{t\to \d^-}\frac{(f')^2}{4} 
\sum_{i=1}\frac{k_i}{h_i - f k_i}  - m_1) = m_1 + 
m_2\tag7.6$$
Note that (7.3), (7.5) and (7.6) determine completely the endpoints
$Z_1$
and $Z_2 = Z_1 - f_\d Z$ which are
$$Z_1 = m_1 Z + Z^\kappa\qquad
Z_2 =  - m_2 Z + Z^\kappa\ . \tag7.7$$
\bigskip
We may now prove the last main theorem. \par
\bigskip
\demo{Proof of Theorem 1.3} The necessity of
(1.3) follows immediately from (7.6) and 
(7.7). \par
To prove (1.4), let us find
 an admissible parametrization $f$ which solves (7.4) 
 with $D=-m_{1}$.
Let us introduce the notation 
$$\Cal F(f) \= \frac{1}{4} 
\sum_i\frac{k_i}{h_i - f k_i}\qquad \Cal H(f) \= f - m_{1}$$
and consider the function $p$ which verifies the differential 
equation $p(f) = f'$. Then (7.4) implies that $p$ verifies
$p' p - p^2 \Cal F(f) + \Cal H(f) = 0$.
So, if we set $p^2 = u$, then 
$$\frac{1}{2} u' - u \Cal F(f) + \Cal H(f) = 0\tag7.8$$
The general solution of (7.8) is
$$u = - 2 e^{2\int  \Cal F(v) dv} 
\left[ \int_0^f \Cal H(v) e^{-2\int {\Cal F(s)}ds}dv + C
\right] = $$
$$ = - 2\frac{1}{\sqrt{\prod_i|h_i - k_i f|}}
\left[ \int_0^f \sqrt{\prod_i|h_i - k_i v|} (v - m_{1}) dv + C
\right]
\tag7.9$$
where the constant $C$ must be determined by the initial
conditions. Since we want that $\lim_{t\to 0^+}f = 0$,
we may immediately set $C=0$.\par
Moreover, since the solution of our
original problem must verify $f(\d) = 
f_\d = m_1+m_2$, $f'(t) >0$ for $t\in ]0, \d[$ and
$f'(0) = f'(\d) = 0$, the solutions
of (7.9) we are interested in verify the conditions
$$u(f) >0 \quad f\in ]0, f_\d[\ ,\qquad
\lim_{f\to 0^+} u(f)= 0 = \lim_{f\to f_\d^-}u(f)$$
This implies that  $f_\d= m_1+m_2$
is the first value after the $0$ for which the following integral vanishes
$$\int_0^{m_1+m_2} \sqrt{\prod_i|h_i - k_i v|}
 (v - m_1) dv = 0 \tag7.10$$
Note that if this occurs, the only values for $f$
on which $u$ vanishes are exactly $0$ and $m_1 + m_2$.\par 
Changing the variable $y=v-m_1$, (7.10) reduces to
$$\int_{-m_1}^{m_2} y\sqrt{\prod_i|(h_i+m_1k_i) - k_i y|}
 dy = 0.\tag 7.11$$
Now, consider the root system defined in \S 1 and take as orthonormal 
basis $\{v_i\}$ the one formed by
the vectors
${1\over 2}(E_\alpha + E_{-\alpha})$ and ${1\over{2i}}(E_\alpha - E_{-\alpha})$.
Then (7.7) implies immediately the necessity of
(1.4).\par
Now, recalling that $u=p^2=(f')^2$, from (7.9) we obtain that
$$t(f) = \int_{0}^f \frac{\root4\of{\prod_i|h_i - k_i s|}}
{\sqrt{
 - 2
 \int_0^s \sqrt{\prod_i |h_i - k_i v|}
 (v - m_1) dv
}
}ds\ .\tag7.12$$
Then set $\d \= t(f_{\d}) = t((m_1+m_2))$. The restriction of the desired
function $f(t)$  on the interval $]0, \d[$ is the inverse function of
(7.12). We now should verify that the function $f$ extends to a 
$C^\infty$-function on $\Bbb R$, which is invariant by the reflections 
at $0$ and $\d$.
It is not difficult, by direct computation, to show that $f$ 
extends to a $C^3$
function which is invariant by the reflections at $0$ and $\d$. 
Now we note that the function $a(t)=f'(t)$ gives rise to a 
$C^2$ Einstein metric $g$; by a 
result of DeTurck and Kazdan (see [Be]), the metric $g$ is real analytic in 
geodesic normal
coordinates, so that the function $f$ is $C^\infty$. \par
The sufficiency of conditions (1), (2) and (3) is clear.\qed
\enddemo

\enddocument

\Refs  
\widestnumber\key{AAAA}  
  

\ref  
\key AA  
\by A.V. Alekseevsky, D.V.Alekseevsky  
\paper G-manifolds with one dimensional orbit space  
\jour Adv. in Sov. Math.  
\vol 8  
\yr 1992   
\pages 1--31  
\endref  
\ref  
\key AA1  
\bysame   
\paper Riemannian G-manifolds with one dimensional orbit space  
\jour Ann. Glob. Anal. and Geom.  
\vol 11  
\yr 1993  
\pages 197--211  
\endref  
\ref
\key AS
\by D.V. Alekseevsky, A. Spiro
\paper Invariant CR-structures on compact homogeneous manifolds
\publ preprint 
\yr 1997
\endref
\ref
\key Be 
\by A.L. Besse
\book Einstein manifolds
\publ Springer-Verlag 
\yr 1986
\endref
\ref  
\key Br  
\by G.E. Bredon  
\book Introduction to compact transformation groups  
\publ Acad. Press N.Y. London  
\yr 1972  
\endref
\ref
\key BFR
\by M. Bordemann, M. Forger, H. R\"omer
\paper Paving the way towards supersymmetric sygma-models
\jour Comm. Math. Physics
\yr 1986
\pages 103--145
\endref
\ref 
\key CN
\by  B.Y. Chen, T. Nagano 
\paper Totally geodesic submanifolds of symmetric spaces, I
\jour Duke Math. J.
\vol 44
\yr 1977
\pages 745--755
\endref
\ref
\key DW
\by A. Dancer, M. Wang
\paper K\"ahler-Einstein metrics of cohomogeneity one and 
bundle constructions for Einstein Hermitia metrics
\publ preprint
\yr 1997
\endref
\ref 
\key GS
\by V. Guillemin, S. Sternberg
\book Symplectic techniques in physics
\publ Cambridge Univ. Press
\yr 1984
\endref
\ref
\key HL
\by W.Y. Hsiang, H.B. Lawson
\paper Minimal submanifolds of low cohomogeneity
\jour J.Diff.Geom.
\vol 5
\yr 1971
\pages 1--38
\endref
\ref
\key HS
\by A. Huckleberry, D. Snow
\paper Almost-homogeneous K\"ahler manifolds with hypersurface 
orbits
\jour Osaka J. math.
\vol 19
\yr 1982
\pages 763--786
\endref
\ref
\key KS
\by N. Koiso, Y. Sakane
\paper Non-homogeneous K\"ahler-Einstein metrics on compact complex manifolds II
\jour Osaka J. math.
\vol 25 
\yr 1988
\pages 933--959
\endref 
\ref
\key PS
\by F.Podest\`a, A. Spiro
\paper K\"ahler manifolds with large isometry group
\jour preprint 
\yr 1997
\endref
\ref
\key Sa
\by Y. Sakane
\paper Examples of compact K\"ahler Einstein manifolds with positive 
Ricci curvature 
\jour Osaka Math. J.
\vol 23
\yr 1986
\pages 585--616
\endref
\ref
\key Ve
\by L. Verdiani
\paper Invariant metrics on cohomogeneity one manifolds
\jour preprint
\yr 1996
\endref
\endRefs
\bye